\documentclass[11pt]{article}
\usepackage{IEEEtrantools}
\usepackage{mathtools, nccmath}
\usepackage{latexsym}
\usepackage{amsfonts}
\usepackage{amssymb}
\usepackage{amsmath,amsfonts,amsthm}
\usepackage{mathrsfs}
\usepackage{enumerate}
\usepackage{dsfont}
\usepackage{hyperref}
\usepackage{graphicx}
\usepackage{epstopdf}
\usepackage{tikz}
\usetikzlibrary{shapes,arrows.meta,positioning}
\usepackage{cleveref}
\usepackage{tikz-cd}
\usepackage{psfrag}

\usepackage{multirow}
\usepackage{tikz}
\usetikzlibrary{arrows,decorations.markings}
\usepackage{bookmark}
\usepackage{hyperref}
\hypersetup{
	colorlinks   = true,
	citecolor    = blue
}

\newcommand{\be}{\begin{equation}}
	\newcommand{\ee}{\end{equation}}
\newcommand{\bea}{\begin{eqnarray}}
	\newcommand{\eea}{\end{eqnarray}}
\newcommand{\bean}{\begin{eqnarray*}}
	\newcommand{\eean}{\end{eqnarray*}}
\newcommand{\brray}{\begin{array}}
	\newcommand{\erray}{\end{array}}
\newcommand{\biearray}{\begin{IEEEarray}{rCl}}
	\newcommand{\eiearray}{\end{IEEEarray}}

\newcommand{\newsection}[1]{\setcounter{equation}{0}
	\setcounter{dfn}{0}
	\section{#1}}

\newtheorem{dfn}{Definition}[section]
\newtheorem{thm}[dfn]{Theorem}
\newtheorem{lmma}[dfn]{Lemma}
\newtheorem{ppsn}[dfn]{Proposition}
\newtheorem{crlre}[dfn]{Corollary}
\newtheorem{xmpl}[dfn]{Example}
\newtheorem{rmrk}[dfn]{Remark}

\newcommand{\bdfn}{\begin{dfn}\rm}
	\newcommand{\bthm}{\begin{thm}}
		\newcommand{\blmma}{\begin{lmma}}
			\newcommand{\bppsn}{\begin{ppsn}}
				\newcommand{\bcrlre}{\begin{crlre}}
					\newcommand{\bxmpl}{\begin{xmpl}}
						\newcommand{\brmrk}{\begin{rmrk}\rm}
							
							\newcommand{\edfn}{\end{dfn}}
						\newcommand{\ethm}{\end{thm}}
					\newcommand{\elmma}{\end{lmma}}
				\newcommand{\eppsn}{\end{ppsn}}
			\newcommand{\ecrlre}{\end{crlre}}
		\newcommand{\exmpl}{\end{xmpl}}
	\newcommand{\ermrk}{\end{rmrk}}

\newcommand{\bbc}{\mathbb{C}}
\newcommand{\bbz}{\mathbb{Z}}
\newcommand{\bbn}{\mathbb{N}}
\newcommand{\bbbt}{\mathbb{T}}
\newcommand{\scrt}{\mathscr{T}}

\newcommand{\cla}{\mathcal{A}}

\newcommand{\clb}{\mathcal{B}}
\newcommand{\clh}{\mathcal{H}}

\def \bbt {\mbox{\boldmath $t$}}

\newcommand{\prf}{\noindent{\it Proof\/}: }
\newcommand{\ots}{\otimes}

\def \qed { \mbox{}\hfill
	$\Box$\vspace{1ex}}

\newcommand{\bq}{\overline{q}}

\title{Sections and Chapters}
\begin{document}

\tikzset{->-/.style={decoration={
  markings,
  mark=at position #1 with {\arrow{>}}},postaction={decorate}}}
  \tikzset{-<-/.style={decoration={
  markings,
  mark=at position #1 with {\arrow{<}}},postaction={decorate}}}

\author{\sc{Debabrata Jana}}
\title{On some algebraic and geometric aspects of the quantum unitary group}
\maketitle


\begin{abstract}
Consider the compact quantum group $U_q(2)$, where $q$ is a non-zero complex deformation parameter such that $|q|\neq 1$. Let $C(U_q(2))$ denote the underlying $C^*$-algebra of the compact quantum group $U_q(2)$. We prove that if $q$ is a non-real complex number and $q^\prime$ is real, then the underlying $C^*$-algebras $C(U_q(2))$ and $C(U_{q^\prime}(2))$ are non-isomorphic. This is in sharp contrast with the case of braided $SU_q(2)$, introduced earlier by Woronowicz et al., where $q$ is a non-zero complex deformation parameter. In another direction, on a geometric aspect of $U_q(2)$, we introduce torus action on the $C^*$-algebra $C(U_q(2))$ and obtain a $C^*$-dynamical system $(C(U_q(2)),\mathbb{T}^3,\alpha)$. We construct a $\mathbb{T}^3$-equivariant spectral triple for $U_q(2)$ that is even and $3^+$-summable. It is shown that the Dirac operator is $K$-homologically nontrivial.	
\end{abstract}
\bigskip

{\bf AMS Subject Classification No.:} 58B32, 58B34, 46L89, 16T20

\hypersetup{linkcolor=blue}
\bigskip

	
\newsection{Introduction}\label{Sec1}

In the theory of compact quantum group (CQG) of Woronowicz \cite{Wor1, Wor2}, the first non-trivial and the most studied example is the $SU_q(2)$ for $q\in\mathbb{R}\setminus\{0\}$. It is widely investigated in the literature through different perspective. In \cite{KasMeyRoyWor-2016aa}, Woronowicz et al. defined a family of $q$-deformations of $SU(2)$ for $q\in\mathbb{C}\setminus\{0\}$. This agrees with the compact quantum group $SU_q(2)$ when $q$ is real but for $q\in\mathbb{C}\setminus\mathbb{R},\,SU_q(2)$ is not a CQG, rather a braided quantum group in a suitable tensor category. In \cite{MeyRoyWor-2016aa}, it is shown that the quantum analogue of the semidirect product construction for groups turns the braided quantum group $SU_q(2)$ into a genuine CQG. This CQG is the coopposite of the compact quantum group $U_q(2)$ defined in \cite{ZhaZha-2005aa}, which we refer as the quantum unitary group. In the case of complex deformation, $U_q(2)$ is `in some sense' the first concrete example of a CQG, since $SU_q(2)$ is no longer a CQG when $q$ is complex.

In this article, we are interested in $U_q(2)$ when $|q|\neq 1$. Note that the $U_q(2)$ is non-Kac when $|q|\neq 1$, whereas it is of Kac-type when $|q|=1$ \cite{ZhaZha-2005aa}. The representation theory of $U_q(2)$ is studied in \cite{GS-20211aa, GS4}, whereas noncommutative geometric aspects have been investigated in \cite{GS-2022aa, GS-20212aa}. It is known that $U_q(2)\cong U_{q'}(2)$ as compact quantum groups if and only if $q'\in\{q,1/q,\overline{q},1/\overline{q}\}$ (Theorem $7.1$, \cite{GS-20211aa}). In this article, our focus is on certain algebraic and geometric aspects of the CQG $U_q(2)$ when $|q|\neq 1$. For $q\neq 0$, let $\theta=\frac{1}{\pi}\arg(q)$. We work under the hypothesis that $\theta$ is irrational. This is along the line of \cite{GS-2022aa, GS-20212aa}, and justification of this is provided in \Cref{Sec2}. Throughout the article, we reserve the notation $C(U_q(2))$ to denote the underlying $C^*$-algebra of the CQG $U_q(2)$.

In the first part of the article, on some algebraic aspects of $U_q(2)$, we investigate the center of $C(U_q(2))$ and observe a striking difference with its $SU$-counterpart. It turns out that if $q$ is a real deformation parameter, then center of $C(U_q(2))$ is $C(\mathbb{T})$, the algebra of continuous functions on torus. However, quite interestingly, when $q$ is a non-real complex deformation parameter, the center of $C(U_q(2))$ (consequently, that of $\mathcal{O}(U_q(2))$) becomes trivial. This helps us to immediately conclude that if $q$ is non-real complex and $q'$ is real, then $C(U_q(2))$ and $C(U_{q'}(2))$ are non-isomorphic $C^*$-algebras. Recall that in the case of braided quantum $SU_q(2)$ introduced in \cite{KasMeyRoyWor-2016aa}, the $C^*$-algebras $C(SU_q(2))$ for different $q$ with $|q|\neq 0,1$ are isomorphic. Therefore, the deformation parameter of the CQG $U_q(2)$ is crucial even at the level of underlying $C^*$-algebras. In another direction, recall that motivated by Jones index for subfactors \cite{J}, Watatani \cite{Wat-1990} had introduced an index for unital inclusion of $C^*$-algebras equipped with a conditional expectation having a quasi-basis, which is a generalization of the Pimsner-Popa basis defined in \cite{Pim-Pop-1986aa}. We consider the homogeneous spaces $U_q(2)/_{\phi}\mathbb{T}$ and $U_q(2)/_{\psi}\mathbb{T}$ of $U_q(2)$ introduced in \cite{GS-2022aa}, and show that both the conditional expectations $E_\phi:C(U_q(2))\to C(U_q(2)/_{\phi}\mathbb{T})$ and $E_\psi:C(U_q(2))\to C(U_q(2)/_{\psi}\mathbb{T})$ have infinite Watatani index.
\smallskip

In the second part of the article, we investigate some geometric aspects of $U_q(2)$. CQGs are rich source of examples in the noncommutative geometry program of Connes \cite{Con1, Con2}. Central notion in Connes’ formulation of Noncommutative geometry is the notion of ``{\em spectral triple}". The first bridge between Woronowicz's CQG and Connes' noncommutative geometry is built through the construction of equivariant $K$-homologically non-trivial spectral triples on $SU_q(2)$ \cite{CP, DLSSV}. In the case of the compact quantum group $U_q(2)$, the representation theory is studied in detail in \cite{GS-20211aa, GS4}, and $K$-homologically non-trivial spectral triple on it that is equivariant under its own comultiplication action has been obtained in \cite{GS-20212aa}, along with the $K$-theory of its underlying $C^*$-algebra. Associated homogeneous spaces have been investigated in \cite{GS-2022aa}. In this article, we obtain a strongly continuous action $\alpha$ of the compact Lie group $\mathbb{T}^3$ on $U_q(2)$, and form the $C^*$-dynamical system $(C(U_q(2)),\mathbb{T}^3,\alpha)$. Let $\clh$ be the separable Hilbert space $\ell^2(\bbn)\otimes\ell^2(\bbz)\otimes\ell^2(\bbz)$. A covariant representation of the $C^*$-dynamical system $(C(U_q(2)),\mathbb{T}^3,\alpha)$ on $\mathcal{H}$ is exhibited, and a $\mathbb{T}^3$-equivariant $3^+$-summable even spectral triple on the Hilbert space $\mathcal{H}\otimes\mathbb{C}^2$ is constructed. Finally, through the index pairing, it is shown that the spectral triple is $K$-homologically non-trivial. Along the process, a relationship with the $\mathbb{T}^2$-equivariant Dirac operator on $SU_q(2)$ obtained in \cite{CP1} is highlighted.
\medskip

Throughout the article, `$i$' will always denote an index and whenever the complex indeterminate appears, we explicitly write it as $\sqrt{-1}$.
	
	
\newsection{Preliminaries}\label{Sec2}
	
We begin by recalling the compact quantum group $U_q(2)$ from \cite{ZhaZha-2005aa}. Let $q$ be a nonzero complex number. The $C^*$-algebra $C(U_q(2))$ is the universal $C^*$-algebra generated by $a,b,D$ satisfying the following relations~:
\begin{IEEEeqnarray}{lCl}\label{relations}
ba&=& q ab, \qquad a^*b=qba^*, \qquad   \qquad \qquad bb^*=b^*b, \qquad \qquad aa^*+bb^*=1, \nonumber \\
aD&=&Da, \qquad bD=q^2|q|^{-2}Db,  \qquad DD^*=D^*D=1,  \qquad a^*a+|q|^2b^*b=1.  
\end{IEEEeqnarray}
The compact quantum group structure is given by the comultiplication $\Delta:C(U_q(2))\longrightarrow C(U_q(2))\otimes C(U_q(2))$ defined as follows~:
\begin{IEEEeqnarray}{lCl} \label{comul}
\Delta(a)=a \otimes a-\bar{q}b \otimes Db^*\quad,\quad\Delta(b)=a\otimes b+b\otimes Da^*\quad,\quad\Delta(D)=D \otimes D. 
\end{IEEEeqnarray}
Let $\mathcal{O}(U_q(2))$ be the $\star$-subalgebra of $C(U_q(2))$ generated by $a,b$ and $D$. The Hopf $\star$-algebra structure on this is given by the following~: 
\begin{IEEEeqnarray*}{lCl}
\mbox{antipode:} \quad S(a)=a^*,\,\,S(b)=-qbD^*,\,\,S(D)=D^*,\,\,S(a^*)=a,\,\,S(b^*)=-(\bar{q})^{-1}b^*D\,,\\ 
\mbox{counit:} \qquad\epsilon(a)=1,\,\,\epsilon(b)=0,\,\,\epsilon(D)=1\,.
\end{IEEEeqnarray*}
$U_q(2)$ is a coamenable compact quantum group that is neither commutative nor cocommutative. For $n,l \in \bbz$ and $ m,k \in \bbn$, define
\begin{IEEEeqnarray}{rCl}\label{basis}
 \langle n,m,k,l \rangle=\begin{cases}
                            a^n b^m(b^*)^k D^l & \mbox{ if } n \geq 0, \cr
                            (a^*)^{-n} b^m(b^*)^k D^l & \mbox{ if } n \leq 0. \cr
                           \end{cases}
\end{IEEEeqnarray} 
\bthm[\cite{ZhaZha-2005aa}]
The set $\{\langle n,m,k,l \rangle: n,l \in \bbz, m,k \in \bbn\}$ forms a linear basis of $\mathcal{O}(U_q(2))$ for all $q\in\mathbb{C}\setminus\{0\}$.
\ethm
\smallskip

In this article, we restrict our attention to $U_q(2)$ when $\,|q|\neq 1$. This is along the line of \cite{GS-20211aa, GS-2022aa, GS-20212aa}.
\begin{ppsn}[Theorem $7.1$, \cite{GS-20211aa}]\label{bs}
The compact quantum groups $U_q(2)$ and $U_{q'}(2)$ are isomorphic as CQG if and only if $q'\in\{q,1/q,\overline{q},1/\overline{q}\}$.
\end{ppsn}
For the purpose of this article, the above classification result allows us to assume that $q\in\bbc$ and $|q|<1$ without loss of generality, and henceforth we do so. Fix any $q\in\mathbb{C}$ such that $0<|q|<1$ and let $\theta=\frac{1}{\pi}\arg{(q)}$.

Consider the Hilbert space $\clh:=\ell^2(\bbn)\otimes\ell^2(\bbz)\otimes\ell^2(\bbz)$. We have the following faithful representation (Proposition $2.1$ in \cite{GS-20211aa}) of $C(U_q(2))$ on $\clh$~:
\begin{IEEEeqnarray}{rCl}\label{representation}
\pi(a)= \sqrt{1-|q|^{2N}}\,V\otimes 1\otimes 1\,,\quad 
\pi(b) =q^N \otimes U\otimes 1\,,\quad
\pi(D) = 1\otimes e^{-2\pi\theta\sqrt{-1}N}\otimes U\,,
\end{IEEEeqnarray}
where $V:e_n\mapsto e_{n+1}$ is the right shift operator acting on $\ell^2(\bbn)$ and $U:e_n\mapsto e_{n+1}$ is the unitary shift operator acting on $\ell^2(\bbz)$. The Haar state $h:C(U_q(2))\longrightarrow\bbc$ is given by the following formula
\begin{IEEEeqnarray}{rCl}\label{haarstate}
h(x)= (1-|q|^2)\sum_{i=0}^{\infty} |q|^{2i} \langle e_{i,0,0}\,,\,\pi(x)e_{i,0,0}\rangle\,,
\end{IEEEeqnarray}
where $\{e_{i,j,k}:i\in\bbn,j,k\in\bbz\}$ denotes the standard orthonormal basis of $\ell^2(\bbn)\otimes\ell^2(\bbz)\otimes\ell^2(\bbz)$. It is shown in (Theorem $2.8$, \cite{GS-20211aa}) that the Haar state is faithful. Note that in \cite{MeyRoyWor-2016aa}, it is shown that at the level of $C^*$-algebras, we have $C(U_q(2))=C(SU_q(2))\boxtimes C(\mathbb{T})$ where $\boxtimes$ is some twisted tensor product, and if $q$ is a real deformation parameter, then $\boxtimes$ coincides with the usual tensor product $\otimes$. Henceforth, we work with the following minor hypothesis (see Section $3$ in \cite{GS-20212aa}).
\medskip

\noindent\textbf{Hypothesis\,:} $\theta$ is irrational where $\theta=\frac{1}{\pi}\arg(q)$.
\medskip

There are two reasons for this additional minor hypothesis. Firstly, in \Cref{Sec3}, we shall be crucially using the irrationality of $\theta$ to compute the center of $C(U_q(2))$ when $q$ is a non-real complex deformation parameter, and secondly, in \Cref{Sec5} we shall crucially use it to show non-triviality of the Diract operator that we construct. The crux here is the simpleness of the noncommutative torus $\mathbb{A}_\theta$ for $\theta$ irrational.
\smallskip

Fix the orthonomal basis $\{e_n:n=0,1,\ldots\}$ of $\ell^2(\mathbb{N})$ recall the bra-ket notation $|e_m\rangle\langle e_n|$ that denotes the rank one projection $e_k\longmapsto e_m\langle e_n,e_k\rangle$. Consider the following operators acting on $\clh$,
\begin{IEEEeqnarray*}{rCl}
a_0=V\otimes 1\otimes 1\quad,\quad b_0=p\otimes U\otimes 1\quad,\quad D_{\theta}=1\otimes e^{-2\sqrt{-1}\pi\theta N} \otimes U\,\,.
\end{IEEEeqnarray*}
Let $C(U_{0,\theta})$ be the $C^*$-subalgebra of $\mathcal{B}(\clh)$ generated by $a_0$, $b_0$ and $D_{\theta}$.
\bppsn[Proposition $3.1$, \cite{GS-20212aa}]
One has $C(U_q(2))=C(U_{0,\theta})$ as $C^*$-algebras.
\eppsn
Let $\scrt:=C^*(V)$ be the Toeplitz algebra. We have the well-known short exact sequence
\[
0 \longrightarrow\mathcal{K}(\ell^2(\bbn))\stackrel{\iota}{\longrightarrow}\scrt\stackrel{\sigma}{\longrightarrow} C(\mathbb{T}) \longrightarrow 0
\]
where $\sigma:V\longmapsto\mathbf{z}$ (here $\mathbf{z}$ denotes the standard unitary generator for $C(\mathbb{T})$). Consider the homomorphism $\tau:C(U_q(2))\longrightarrow C(\bbbt)\otimes\clb(\ell^2(\bbz)\otimes\ell^2(\bbz))$ given by $\tau=\sigma\otimes 1\otimes 1$, and let
\begin{IEEEeqnarray*}{rCl}
\mathcal{I}_\theta &=& \mbox{ the closed two-sided ideal of } C(U_q(2)) \mbox{ generated by } b_0 \mbox{ and } b_0^*\,,\\
\mathcal{B}_\theta &=& C^*\big(\{\tau(a_0)\,,\,\tau(D_\theta)\}\big)=C^*\big(\{\mathbf{z}\otimes 1\otimes 1\,,\,1\otimes e^{-2\sqrt{-1}\pi\theta N}\otimes U\}\big)\,.
\end{IEEEeqnarray*}

\bppsn[Proposition $3.2$, \cite{GS-20212aa}]\label{short exact}
The following chain of $C^*$-algebras
\[
0 \longrightarrow\mathcal{I}_\theta\stackrel{\iota}{\longrightarrow} C(U_q(2))\stackrel{\tau}{\longrightarrow}\mathcal{B}_\theta\longrightarrow 0
\]
is a short exact sequence, where `$\iota$' denotes the inclusion map.
\eppsn

Let $\mathcal{C}_\theta$ be the $C^*$-subalgebra of $\mathcal{B}(\ell^2(\bbz)\otimes\ell^2(\bbz))$ generated by $U\otimes 1$ and $e^{-2\sqrt{-1}\pi\theta N}\otimes U$. Since $\theta$ is irrational, by the universality and simpleness of the noncommutative torus $\mathbb{A}_\theta$, we get that $\mathcal{C}_\theta\cong\mathbb{A}_\theta$ as $C^*$-algebras. Moreover, the following holds.
\blmma[Lemma $3.3$, \cite{GS-20212aa}]\label{kernel is NC torus}
$\mathcal{I}_\theta=\mathcal{K}(\ell^2(\bbn))\otimes\mathcal{C}_\theta\,$.
\elmma


\newsection{On isomorphism classes of the $C^*$-algebras $C(U_q(2))$}\label{Sec3}

Let $C(U_q(2))$ be the underlying $C^*$-algebra of the compact quantum group $U_q(2)$ and $\mathcal{O}(U_q(2))$ be the Hopf $\star$-algebra generated by $a,b,D$. In the case of braided quantum group $SU_q(2),\,q\in\bbc\setminus\{0\},$ introduced in \cite{KasMeyRoyWor-2016aa}, the underlying $C^*$-algebras $C(SU_q(2))$ for different $q$ with $|q|\neq 0,1$ are isomorphic (see Theorem $2.3$ in \cite{KasMeyRoyWor-2016aa}). In this section, we show that this does not hold in the case of $U_q(2)$. More precisely, we prove the following result.

\begin{thm}\label{main thm}
Let $q$ be a non-real complex number such that $|q|\neq 1$ and $q'\in\mathbb{R}$ be a non-zero real number. Let $\theta=\frac{1}{\pi}\arg{(q)}$ and $\theta$ is irrational. Then, the $C^*$-algebras $C(U_q(2))$ and $C(U_{q'}(2))$ are non-isomorphic.
\end{thm}

We prove the above theorem by investigating the center of the $C^*$-algebras $C(U_q(2))$ for different parameter values of $q$. For any algebra $\mathscr{A}$, we denote its center by $\mathcal{Z}(\mathscr{A})$.

\bthm\label{cuq2 center}
When $q$ is a real deformation parameter, the center $\,\mathcal{Z}(C(U_q(2)))$ is isomorphic to $C(\mathbb{T})$, and when $q$ is a non-real complex deformation parameter, $\,\mathcal{Z}(C(U_q(2)))=\bbc$.
\ethm
\prf Recall from \Cref{bs} that $U_q(2)$ and $U_{q'}(2)$ are isomorphic as compact quantum groups if and only if $q'\in\{q,1/q,\overline{q},1/\overline{q}\}$. If two compact quantum groups are isomorphic as CQG, then it is necessary that their underlying $C^*$-algebras are isomorphic. Therefore, without loss of generality, we can assume that $q\in\bbc$ and $0<|q|<1$, that is, $q$ lies in the open unit disk.

When $q$ is a non-zero real number in $(-1,1)$, we have $C(U_q(2))=C(SU_q(2))\otimes C(\mathbb{T})$ as $C^*$-algebras (see Theorem $2.1$ in \cite{Z}). It is well-known that $C(SU_q(2))$ has trivial center, and thus by \cite{Haydon-1973aa} we have $\mathcal{Z}\big(C(U_q(2))\big)=C(\mathbb{T})$.

Now, suppose that $q$ is a non-real complex deformation parameter such that $|q|<1$. Consider the Hilbert space $\mathcal{H}=\ell^2(\mathbb{N})\otimes\ell^2(\mathbb{Z})\otimes\ell^2(\mathbb{Z})$ and its standard orthonormal basis $\{e_{d,i,j}:d\in\bbn; i,j\in\bbz\}$. Recall the faithful representation $\pi$ of $C(U_q(2))$ on $\mathcal{H}$ from \Cref{representation}. Motivated by the analysis carried out in \cite{Soltan-2016aa}, we first let $T\in C(U_q(2))^\prime\subseteq\mathcal{B}(\mathcal{H})$, where $C(U_q(2))^\prime$ denotes the commutant in $\mathcal{B}(\mathcal{H})$. The general expression of $T$ is the following
\[T(e_{d,i,j})=\sum_{d^\prime,i^\prime,j^\prime}\alpha_{d^\prime,i^\prime,j^\prime}^{d,i,j}\,e_{d^\prime,i^\prime,j^\prime}\]
where each $\alpha_{d^\prime,i^\prime,j^\prime}^{d,i,j}\in\bbc$. Fix $(d,i,j)\in\mathbb{N}\times\mathbb{Z}\times\mathbb{Z}$. Omitting the representation symbol $\pi$ for brevity, using \Cref{representation}, we observe that the following holds\,:
\begin{IEEEeqnarray*}{rCl}
Tb(e_{d,i,j}) &=& \sum_{d^\prime,i^\prime,j^\prime}q^d\alpha_{d^\prime,i^\prime,j^\prime}^{d,i+1,j}\,e_{d^\prime,i^\prime,j^\prime}=\sum_{d^\prime,i^\prime,j^\prime}q^d\alpha_{d^\prime,i^\prime+1,j^\prime}^{d,i+1,j}\,e_{d^\prime,i^\prime+1,j^\prime}\,,\\
bT(e_{d,i,j}) &=& \sum_{d^\prime,i^\prime,j^\prime}q^{d^\prime}\alpha_{d^\prime,i^\prime,j^\prime}^{d,i,j}\,e_{d^\prime,i^\prime+1,j^\prime}\,,\\
Tb^*(e_{d,i,j}) &=& \sum_{d^\prime,i^\prime,j^\prime}(\overline{q})^d\alpha_{d^\prime,i^\prime,j^\prime}^{d,i-1,j}\,e_{d^\prime,i^\prime,j^\prime}=\sum_{d^\prime,i^\prime,j^\prime}(\overline{q})^d\alpha_{d^\prime,i^\prime-1,j^\prime}^{d,i-1,j}\,e_{d^\prime,i^\prime-1,j^\prime}\,,\\
b^*T(e_{d,i,j}) &=& \sum_{d^\prime,i^\prime,j^\prime}(\overline{q})^{d^\prime}\alpha_{d^\prime,i^\prime,j^\prime}^{d,i,j}\,e_{d^\prime,i^\prime-1,j^\prime}\,,\\
Ta(e_{d,i,j}) &=& \sum_{d^\prime,i^\prime,j^\prime}\sqrt{1-|q|^{2d+2}}\,\alpha_{d^\prime,i^\prime,j^\prime}^{d+1,i,j}\,e_{d^\prime,i^\prime,j^\prime}=\sum_{d^\prime,i^\prime,j^\prime}\sqrt{1-|q|^{2d+2}}\,\alpha_{d^\prime+1,i^\prime,j^\prime}^{d+1,i,j}\,e_{d^\prime+1,i^\prime,j^\prime}\,,\\
aT(e_{d,i,j}) &=& \sum_{d^\prime,i^\prime,j^\prime}\sqrt{1-|q|^{2d^\prime+2}}\,\alpha_{d^\prime,i^\prime,j^\prime}^{d,i,j}\,e_{d^\prime+1,i^\prime,j^\prime}\,,\\
Ta^*(e_{d,i,j}) &=& \sum_{d^\prime,i^\prime,j^\prime}\sqrt{1-|q|^{2d}}\,\alpha_{d^\prime,i^\prime,j^\prime}^{d-1,i,j}\,e_{d^\prime,i^\prime,j^\prime}=\sum_{d^\prime,i^\prime,j^\prime}\sqrt{1-|q|^{2d}}\,\alpha_{d^\prime-1,i^\prime,j^\prime}^{d-1,i,j}\,e_{d^\prime-1,i^\prime,j^\prime}\,,\\
a^*T(e_{d,i,j}) &=& \sum_{d^\prime,i^\prime,j^\prime}\sqrt{1-|q|^{2d^\prime}}\,\alpha_{d^\prime,i^\prime,j^\prime}^{d,i,j}\,e_{d^\prime-1,i^\prime,j^\prime}\,,\\
TD(e_{d,i,j}) &=& \sum_{d^\prime,i^\prime,j^\prime}e^{-2\pi \sqrt{-1}i\theta}\,\alpha_{d^\prime,i^\prime,j^\prime}^{d,i,j+1}\,e_{d^\prime,i^\prime,j^\prime}=\sum_{d^\prime,i^\prime,j^\prime}e^{-2\pi \sqrt{-1}i\theta}\,\alpha_{d^\prime,i^\prime,j^\prime+1}^{d,i,j+1}\,e_{d^\prime,i^\prime,j^\prime+1}\,,\\
DT(e_{d,i,j}) &=& \sum_{d^\prime,i^\prime,j^\prime}e^{-2\pi \sqrt{-1}i^\prime\theta}\,\alpha_{d^\prime,i^\prime,j^\prime}^{d,i,j}\,e_{d^\prime,i^\prime,j^\prime+1}\,,\\
TD^*(e_{d,i,j}) &=& \sum_{d^\prime,i^\prime,j^\prime}e^{2\pi \sqrt{-1}i\theta}\,\alpha_{d^\prime,i^\prime,j^\prime}^{d,i,j-1}\,e_{d^\prime,i^\prime,j^\prime}=\sum_{d^\prime,i^\prime,j^\prime}e^{2\pi \sqrt{-1}i\theta}\,\alpha_{d^\prime,i^\prime,j^\prime-1}^{d,i,j-1}\,e_{d^\prime,i^\prime,j^\prime-1}\,,\\
D^*T(e_{d,i,j}) &=& \sum_{d^\prime,i^\prime,j^\prime}e^{2\pi \sqrt{-1}i^\prime\theta}\,\alpha_{d^\prime,i^\prime,j^\prime}^{d,i,j}\,e_{d^\prime,i^\prime,j^\prime-1}\,.
\end{IEEEeqnarray*}
Thus, we have the following set of equations\,:
\begin{IEEEeqnarray}{rCl}\label{eq}
\alpha_{d^\prime,i^\prime,j^\prime}^{d,i,j} &=& (\overline{q})^{d-d^\prime}\alpha_{d^\prime,i^\prime-1,j^\prime}^{d,i-1,j}\hspace*{2cm}\big(\mbox{as } Tb^*=b^*T\big)\nonumber\\
&=& q^{d-d^\prime}\alpha_{d^\prime,i^\prime+1,j^\prime}^{d,i+1,j}\hspace*{2.3cm}\big(\mbox{as } Tb=bT\big)\nonumber\\
&=& \sqrt{\frac{1-|q|^{2d+2}}{1-|q|^{2d^\prime+2}}}\,\alpha_{d^\prime+1,i^\prime,j^\prime}^{d+1,i,j}\qquad\big(\mbox{as } Ta=aT\big)\nonumber\\
&=& \sqrt{\frac{1-|q|^{2d}}{1-|q|^{2d^\prime}}}\,\alpha_{d^\prime-1,i^\prime,j^\prime}^{d-1,i,j}\hspace*{1.2cm}\big(\mbox{as } Ta^*=a^*T\big)\nonumber\\
&=& e^{2\pi \sqrt{-1}(i^\prime-i)\theta}\,\alpha_{d^\prime,i^\prime,j^\prime+1}^{d,i,j+1}\qquad\big(\mbox{as } TD=DT\big)\nonumber\\
&=& e^{2\pi \sqrt{-1}(i-i^\prime)\theta}\,\alpha_{d^\prime,i^\prime,j^\prime-1}^{d,i,j-1}\qquad\big(\mbox{as } TD^*=D^*T\big)\,.
\end{IEEEeqnarray}
From the first two equations in \Cref{eq} we get the following
\[
\alpha_{d^\prime,i^\prime,j^\prime}^{d,i,j}=q^{d-d^\prime}\alpha_{d^\prime,i^\prime+1,j^\prime}^{d,i+1,j}=|q|^{2(d-d^\prime)}\alpha_{d^\prime,i^\prime,j^\prime}^{d,i,j}\,.
\]
Therefore, $\alpha_{d^\prime,i^\prime,j^\prime}^{d,i,j}=0$ for $d\neq d^\prime$. Hence, the third and the fourth equations in \Cref{eq}, that is, $\alpha_{d^\prime,i^\prime,j^\prime}^{d,i,j}=\sqrt{\frac{1-|q|^{2d+2}}{1-|q|^{2d^\prime+2}}}\,\alpha_{d^\prime+1,i^\prime,j^\prime}^{d+1,i,j}$ and $\alpha_{d^\prime,i^\prime,j^\prime}^{d,i,j}=\sqrt{\frac{1-|q|^{2d}}{1-|q|^{2d^\prime}}}\,\alpha_{d^\prime-1,i^\prime,j^\prime}^{d-1,i,j}$ together imply the following
\begin{IEEEeqnarray}{rCl}\label{eq1}
\alpha_{d+r,i^\prime,j^\prime}^{d+r,i,j} &=& \alpha_{d,i^\prime,j^\prime}^{d,i,j}
\end{IEEEeqnarray}
for all $r\in\bbz$. Moreover, the first two equations $\alpha_{d^\prime,i^\prime,j^\prime}^{d,i,j}=q^{d-d^\prime}\alpha_{d^\prime,i^\prime+1,j^\prime}^{d,i+1,j}$ and $\alpha_{d^\prime,i^\prime,j^\prime}^{d,i,j}=(\overline{q})^{d-d^\prime}\alpha_{d^\prime,i^\prime-1,j^\prime}^{d,i-1,j}$ together imply the following
\begin{IEEEeqnarray}{rCl}\label{eq2}
\alpha_{d,i^\prime+r,j^\prime}^{d,i+r,j} &=& \alpha_{d,i^\prime,j^\prime}^{d,i,j}
\end{IEEEeqnarray}
for all $r\in\bbz$. Therefore, using \Cref{eq1,eq2}, together with the fact that $\alpha_{d^\prime,i^\prime,j^\prime}^{d,i,j}=0$ for $d\neq d^\prime$, we get the following
\begin{IEEEeqnarray}{rCl}\label{a}
\alpha_{d^\prime,i^\prime,j^\prime}^{d,i,j}=\delta_{dd^\prime}\alpha_{0,i^\prime,j^\prime}^{0,i,j} &=& \delta_{dd^\prime}\alpha_{0,0,j^\prime}^{0,i-i^\prime,j}\nonumber\\
&=& \delta_{dd^\prime}\,e^{2\pi \sqrt{-1}j^\prime(i-i^\prime)\theta}\,\alpha_{0,0,0}^{0,i-i^\prime,j-j^\prime}\,,
\end{IEEEeqnarray}
where the last equality in above follows from the last equation $\alpha_{d^\prime,i^\prime,j^\prime}^{d,i,j}=e^{2\pi \sqrt{-1}(i-i^\prime)\theta}\,\alpha_{d^\prime,i^\prime,j^\prime-1}^{d,i,j-1}$ in \Cref{eq}.
	
Now, if we choose $T\in\mathcal{B}(\mathcal{H})$ from the center $\mathcal{Z}(C(U_q(2)))\subseteq C(U_q(2))^\prime$, then using the density of $\mathcal{O}(U_q(2))$ in $C(U_q(2))$, we have
\begin{IEEEeqnarray*}{rCl}
T &=& \lim_{\lambda\rightarrow\infty}\,\sum_{(n,m,k,l)\in\Gamma}C^{\lambda}_{n,m,k,l}\,\langle n,m,k,l\rangle,
\end{IEEEeqnarray*}
where $\langle n,m,k,l\rangle$ is the linear basis of $\mathcal{O}(U_q(2))$ defined in \Cref{basis}, and $\Gamma$ denotes finite subset of $\bbz\times\bbn\times\bbn\times\bbz$. Hence for $\epsilon>0$, there exists $\lambda_0\in\mathbb{N}$ such that for any $\lambda\geq\lambda_0$ we have the following
\begin{IEEEeqnarray*}{rCl}
\Big|\Big|\,T-\sum_{(n,m,k,l)\in\Gamma}C^{\lambda}_{n,m,k,l}\langle n,m,k,l\rangle\,\Big|\Big|<\epsilon/2\,.
\end{IEEEeqnarray*}
Then, using \Cref{a} we get the following
\begin{IEEEeqnarray*}{rCl}
(\epsilon/2)^2 &>& \Big|\Big|\Big(T-\sum_{(n,m,k,l)\in\Gamma}C^{\lambda}_{n,m,k,l}\langle n,m,k,l\rangle\Big)(e_{p,0,0})\,\Big|\Big|^2\\
&=& \Big|\Big|\sum_{t,s\in\mathbb{Z}}\alpha_{0,0,0}^{0,-s,-t}exp(-2\pi \sqrt{-1}st\theta)e_{p,s,t}-\sum_{(n,m,k,l)\in\Gamma}C^{\lambda}_{n,m,k,l}\langle n,m,k,l\rangle(e_{p,0,0})\,\Big|\Big|^2\\
&=& \Big|\Big|\sum_{t,s\in\mathbb{Z}}\alpha_{0,0,0}^{0,-s,-t}exp(-2\pi \sqrt{-1}st\theta)e_{p,s,t}\\
&  & -\sum_{\underset{-p\le n\le 0}{(n,m,k,l)\in\Gamma}}C^{\lambda}_{n,m,k,l}q^{pm}\bq^{pk}\prod_{c\,=0\,}^{-(n+1)}\sqrt{1-|q|^{2(p-c)}}e_{p+n,m-k,l}\\
&  & -\sum_{\underset{1\le n<\infty}{(n,m,k,l)\in\Gamma}}C^{\lambda}_{n,m,k,l}q^{pm}\bq^{pk}\prod_{c\,=\,1}^n\sqrt{1-|q|^{2(p+c)}}e_{p+n,m-k,l}\,\Big|\Big|^2\\
&=& \Big|\Big|\sum_{t,s\in\mathbb{Z}}\alpha_{0,0,0}^{0,-s,-t}exp(-2\pi \sqrt{-1}st\theta)e_{p,s,t}-\sum_{\underset{n=0,\,m-k=s,\,s\in\mathbb{Z}}{(n,m,k,l)\in\Gamma}}C^{\lambda}_{n,m,k,l}q^{pm}\bq^{pk}e_{p,s,l}\\
&  & -\sum_{\underset{-p\le n\le -1\,;\,m-k=s,\,s\in\bbz}{(n,m,k,l)\in\Gamma}}C^{\lambda}_{n,m,k,l}q^{pm}\bq^{pk}\prod_{c\,=\,0}^{-(n+1)}\sqrt{1-|q|^{2(p-c)}}e_{p+n,s,l}\\
&  & -\sum_{\underset{1\le n<\infty\,;\,m-k=s,\,s\in \bbz}{(n,m,k,l)\in\Gamma}}C^{\lambda}_{n,m,k,l}q^{pm}\bq^{pk}\prod_{c\,=\,1}^n\sqrt{1-|q|^{2(p+c)}}e_{p+n,s,l}\,\Big|\Big|^2\\
&=& \Big|\Big|\sum_{t,s\in\mathbb{Z}}\Big(\alpha_{0,0,0}^{0,-s,-t}exp(-2\pi \sqrt{-1}st\theta)-\sum_{\underset{m-k=s}{m,k\in\mathbb{N}}}C^{\lambda}_{0,m,k,t}q^{pm}\bq^{pk}\Big)e_{p,s,t}\\
&  & -\sum_{\underset{-p\le n\le -1\,;\,m-k=s,\,s\in\bbz}{(n,m,k,l)\in\Gamma}}C^{\lambda}_{n,m,k,l}q^{pm}\bq^{pk}\prod_{c\,=\,0}^{-(n+1)}\sqrt{1-|q|^{2(p-c)}}e_{p+n,s,l}\\
&  & -\sum_{\underset{1\le n<\infty\,;\,m-k=s,\,s\in \bbz}{(n,m,k,l)\in\Gamma}}C^{\lambda}_{n,m,k,l}q^{pm}\bq^{pk}\prod_{c\,=\,1}^n\sqrt{1-|q|^{2(p+c)}}e_{p+n,s,l}\,\Big|\Big|^2\\
&=& \sum_{t,s\in\mathbb{Z}}\Big|\Big(\alpha_{0,0,0}^{0,-s,-t}exp(-2\pi \sqrt{-1}st\theta)-\sum_{\underset{m-k=s}{m,k\in\mathbb{N}}}C^{\lambda}_{0,m,k,t}q^{pm}\bq^{pk}\Big)\,\Big|^2\\
&  & +\sum_{\underset{-p\le n\le-1}{n,l,s\in\mathbb{Z}}}\Big|\sum_{\underset{m-k=s}{m,k\in\mathbb{N}}}C^{\lambda}_{n,m,k,l}q^{pm}\bq^{pk}\prod_{c\,=\,0}^{-(n+1)}\sqrt{1-|q|^{2(p-c)}}\,\Big|^2\\
&  & +\sum_{\underset{1\le n<\infty}{n,l,s\in\mathbb{Z}}}\Big|\sum_{\underset{m-k=s}{m,k\in\mathbb{N}}}C^{\lambda}_{n,m,k,l}q^{pm}\bq^{pk}\prod_{c\,=\,1}^n\sqrt{1-|q|^{2(p+c)}}\,\Big|^2\\
&\geq& \sum_{t,s\in\mathbb{Z}}\Big|\alpha_{0,0,0}^{0,-s,-t}exp(-2\pi \sqrt{-1}st\theta)-\sum_{\underset{m-k=s}{m,k\in\mathbb{N}}}C^{\lambda}_{0,m,k,t}q^{pm}\bq^{pk}\Big|^2\\
&=& \sum_{t,s\in\mathbb{Z}}\Big|\alpha_{0,0,0}^{0,-s,-t}-\sum_{\underset{m-k=s}{m,k\in\mathbb{N}}}C^{\lambda}_{0,m,k,t}q^{pm}\bq^{pk}exp(2\pi \sqrt{-1}st\theta)\Big|^2\,.
\end{IEEEeqnarray*}
Therefore, $\forall\,\epsilon>0,\,\forall\,p\in\mathbb{N},\,\forall\,\lambda\in\mathbb{N}$ with $\lambda\geq\lambda_0,\,\forall\,t,s\in\mathbb{Z}$ we get the following
\begin{IEEEeqnarray*}{rCl}
\Big|\alpha_{0,0,0}^{0,-s,-t}-\sum_{\underset{m-k=s}{m,k\in\mathbb{N}}}C^{\lambda}_{0,m,k,t}q^{pm}\bq^{pk}exp(2\pi \sqrt{-1}st\theta)\Big| &<& \epsilon/2\,.
\end{IEEEeqnarray*}
Take $t\in\mathbb{Z},\,s\in\mathbb{Z}\setminus\{0\}$ and $\epsilon>0$. Since,
\begin{IEEEeqnarray*}{rCl}
\lim_{p\rightarrow\infty}\sum_{\underset{m-k=s}{m,k\in\mathbb{N}}}C^{\lambda}_{0,m,k,t}q^{pm}\bq^{pk}exp(2\pi \sqrt{-1}st\theta) &=& 0
\end{IEEEeqnarray*}
as $|q|<1,\,\exists\,p_0\in\mathbb{N}$ such that $\forall\,p\geq p_0$ we have the following,
\begin{IEEEeqnarray*}{rCl}
\Big|\sum_{\underset{m-k=s}{m,k\in\mathbb{N}}}C^{\lambda}_{0,m,k,t}q^{pm}\bq^{pk}exp(2\pi \sqrt{-1}st\theta)\Big| &<& \epsilon/2\,.
\end{IEEEeqnarray*}
Thus, for given $\epsilon>0,\,s\in\mathbb{Z}\setminus\{0\},\,t\in\mathbb{Z}\,\,\exists\,\lambda_0,p_0\in\mathbb{N}$ such that $\forall\,\lambda\geq\lambda_0$ and $p\geq p_0$ we have the following,
\begin{IEEEeqnarray*}{rCl}
\Big|\alpha_{0,0,0}^{0,-s,-t}\Big| &\leq& \Big|\sum_{\underset{m-k=s}{m,k\in\mathbb{N}}}C^{\lambda}_{0,m,k,t}q^{pm}\bq^{pk}exp(2\pi \sqrt{-1}st\theta)-\alpha_{0,0,0}^{0,-s,-t}\Big|\\
&   & +\Big|\sum_{\underset{m-k=s}{m,k\in\mathbb{N}}}C^{\lambda}_{0,m,k,t}q^{pm}\bq^{pk}exp(2\pi \sqrt{-1}st\theta)\Big|\\
&<& \epsilon/2+\epsilon/2=\epsilon\,.
\end{IEEEeqnarray*}
So, any $T\in\mathcal{Z}(C(U_q(2)))\subseteq\mathcal{B}(\mathcal{H})$ is of the following form
\begin{IEEEeqnarray*}{rCl}
T(e_{d,i,j}) &=& \sum_{j^\prime\in\mathbb{Z}}\alpha_{0,0,0}^{0,0,j-j^\prime}\,e_{d,i,j^\prime}\\
&=& \alpha_{0,0,0}^{0,0,0}e_{d,i,j}+\sum_{j^\prime\in\mathbb{Z}}\alpha_{0,0,0}^{0,0,j-j^\prime}(1-\delta_{jj^\prime})\,e_{d,i,j^\prime}\,.
\end{IEEEeqnarray*}
Thus, any operator in $\mathcal{Z}(C(U_q(2)))$ must be of the following form\,:
\begin{IEEEeqnarray}{rCl}\label{operator form1}
T\,:\,e_{d,i,j} &\longmapsto& \sum_{j^\prime\in\mathbb{Z}}\omega^{0,0,j^\prime}_{j^\prime}(T)\big(I\otimes I\otimes (U^*)^{j^\prime}\big)(e_{d,i,j})\\
&=& \Big(I\otimes\Big(\sum_{j^\prime\in\mathbb{Z}}\omega^{0,j^\prime}_{j^\prime}(T)\big(I\otimes (U^*)^{j^\prime}\big)\Big)\Big)(e_{d,i,j})
\end{IEEEeqnarray}
where $U^*$ is the left shift operator acting on $\ell^2(\mathbb{Z})$, as in \Cref{representation}, and
\begin{IEEEeqnarray}{rCl}\label{co-efficient1}
\omega^{0,0,j^\prime}_{j^\prime}(T) &=& \langle e_{0,0,0}\,,\,T(e_{0,0,j^\prime})\rangle\,.
\end{IEEEeqnarray}
It is easy to check that such $T$ is a normal operator. Moreover,
\begin{IEEEeqnarray*}{rCl}
T^*(e_{d,i,j}) &=& \sum_{j^\prime\in\mathbb{Z}}\overline{\omega^{0,0,j^\prime}_{j^\prime}(T)}\,\big(I\otimes I\otimes U^{j^\prime}\big)(e_{d,i,j})\\
&=& \sum_{j^\prime\in\mathbb{Z}}\overline{\omega^{0,0,-j^\prime}_{-j^\prime}(T)}\big(I\otimes I\otimes (U^*)^{j^\prime}\big)(e_{d,i,j})\,.
\end{IEEEeqnarray*}
Therefore, such an operator $T$ is self-adjoint if and only if for each $j\in\bbz$ we have $\,\omega^{0,0,-j}_{-j}(T)=\overline{\omega^{0,0,j}_j(T)}$.
\smallskip
	
Now, by \Cref{representation} it is easy to see that any operator $T=I\otimes T_1$, where $T_1=\sum_{j^\prime\in\mathbb{Z}}\omega^{0,j^\prime}_{j^\prime}(T)\big(I\otimes (U^*)^{j^\prime}\big)$, satisfying \Cref{operator form1,co-efficient1} commutes with $\pi(a),\,\pi(b)$ and $\pi(D)$, that is, $T=I\otimes T_1\in\big(C(U_q(2))\big)^\prime\cap\mathcal{B}(\mathcal{H})$. If $T$ is in the center of $C(U_q(2))$, then by \Cref{kernel is NC torus}, $(p\otimes 1\otimes 1)T=p\otimes T_1$ is in the ideal $\mathcal{I}_\theta\subset C(U_q(2))$, where $p=|e_0\rangle\langle e_0|$, the rank one projection onto $\bbc e_0$. Therefore, we see that $T_1\in\mathcal{C}_\theta$, the unital $C^*$-subalgebra of $\mathcal{B}(\ell^2(\bbz)\otimes\ell^2(\bbz))$ generated by $U\otimes 1$ and $e^{-2\sqrt{-1}\pi\theta N}\otimes U$. However, observe that $T_1\in\mathcal{C}_\theta^\prime\cap\mathcal{B}(\ell^2(\bbz)\otimes\ell^2(\bbz))$ as it commutes with both the generators of $\mathcal{C}_\theta$, and consequently we have $T_1\in\mathcal{C}_\theta\cap\mathcal{C}_\theta^\prime$, that is, $T_1\in\mathcal{Z}(\mathcal{C}_\theta)$. However, $\mathcal{C}_\theta$ is isomorphic to the noncommutative torus $\mathbb{A}_\theta$ and by our hypothesis that $\theta$ is irrational, we get that $\mathcal{C}_\theta$ is a unital simple $C^*$-algebra. Hence, $\mathcal{Z}(\mathcal{C}_\theta)$ must be trivial, and consequently, $T_1$ must be scalar. This finally concludes that any $T\in\mathcal{Z}(C(U_q(2)))$ must be scalar, when $q$ is a non-real complex deformation parameter, which finishes the proof.\qed

\begin{crlre}
The center of $\,\mathcal{O}(U_q(2))$ is trivial when $q$ is a non-real complex deformation parameter and it is non-trivial when $q$ is a nonzero real parameter.
\end{crlre}
\begin{prf}
If $q$ is a non-real complex number, $\mathcal{Z}\big(\mathcal{O}(U_q(2))\big)=\mathbb{C}$ immediately follows from \Cref{cuq2 center}. When $q$ is real, it is clear from \Cref{representation} that $\pi(D)=I\otimes I\otimes U$. Thus, the $\star$-subalgebra generated by $D$ is contained in $\mathcal{Z}\big(\mathcal{O}(U_q(2))\big)$ when $q$ is real.\qed
\end{prf}
\medskip

\noindent\textbf{Proof of Theorem \ref{main thm}\,:} If $q$ is a non-real complex number, then the $C^*$-algebra $C(U_q(2))$ has trivial center by Theorem \ref{cuq2 center}. On the other hand, if $q^\prime\in\mathbb{R}\setminus\{0\}$, then the center of the $C^*$-algebra $C(U_{q^\prime}(2))$ is isomorphic to $C(\mathbb{T})$. Therefore, $C(U_q(2))$ and $C(U_{q^\prime}(2))$ are non-isomorphic.\qed
\smallskip

When $q\in(-1,1)\setminus\{0\}$, we have $C(U_q(2))=C(SU_q(2))\otimes C(\mathbb{T})$ (see Theorem $2.1$ in \cite{Z}), and since the $C^*$-algebras $C(SU_q(2))$ are all isomorphic when $0<|q|<1$, same holds for $C(U_q(2))$. We end this section with the following question\,:
\smallskip

\noindent\textbf{Question\,:} Given $q,q'$ two distinct non-real complex deformation parameters, is $C(U_q(2))\cong C(U_{q'}(2))$ as $C^*$-algebras?
\smallskip


\newsection{Homogeneous spaces of $U_q(2)$ and Watatani index}\label{Sec4}

In this section, we compute the Watatani indices for the homogeneous spaces $C(U_q(2)/_{\phi}\bbbt)$ and $C(U_q(2)/_{\psi}\bbbt)$ of $U_q(2)$ investigated in \cite{GS-2022aa}. We begin by recalling some definitions and results from \cite{Wat-1990} which will be used in this section. Consider a unital inclusion of $C^*$-algebras $\mathcal{B}\subset\mathcal{A}$.
\bdfn(\cite{Wat-1990}).
For a unital inclusion of $C^*$-algebras $\mathcal{B}\subset\mathcal{A}$, a conditional expectation $E:\mathcal{A}\to\mathcal{B}$ is a surjective linear map satisfying the following
\[
E(ba)=bE(a)\quad,\quad E(ab)=E(a)b\quad,\quad E(b)=b
\]
for $b\in\mathcal{B}$ and $a\in\mathcal{A}$. Moreover, we assume that $E$ is positive, that is, $E$ is a projection of norm one.
\edfn 

\bdfn (\cite{Wat-1990}).
For a unital inclusion of $C^*$-algebras $\mathcal{B}\subset\mathcal{A}$, a conditional expectation $E:\mathcal{A}\to\mathcal{B}$ is said to be of index-finite type if there exists a finite set $\{(u_j,v_j):1\leq j\leq n\}\subseteq\mathcal{A}\times\mathcal{A}$ such that
\[
x= \sum_iu_iE(v_ix)=\sum_iE(xu_i)v_i
\]
for any $x \in\mathcal{A}$. Such a set $\{(u_1,v_1),(u_2,v_2),\ldots,(u_n,v_n)\}$ is called a quasi-basis for $E$. The Watatani index of $E$ is defined by 
\[
\mbox{Index}(E):=\sum_iu_i v_i
\]
and is in fact independent of the quasi-basis. Moreover, one can choose a quasi-basis with $v_i=u_i^*$ for all $i$.
\edfn

\bppsn [\cite{Wat-1990}]\label{indexcondition}
For a unital inclusion of $C^*$-algebras $\mathcal{B}\subset\mathcal{A}$, let $E:\mathcal{A}\to\mathcal{B}$ be a conditional expectation. If $E$ is of index-finite type, then $c(E)\geq ||\mathrm{Index } (E)||^{-1}$ where
\[
c(E):=\sup\{t\in\mathbb{R}^+\,:\,E(a)\geq ta\,\,\forall\,a\in\mathcal{A}_+\}\,.
\]
\eppsn

In the case of subfactors, $c(E)$ is called the Pimsner-Popa probabilistic constant \cite{Pim-Pop-1986aa}. Let $\bbt:z\mapsto z$ denote the standard unitary generator of $C(\bbbt)$ and $\Delta_{\bbbt}$ be the standard comultiplication on $C(\bbbt)$. Recall from \cite{GS-2022aa} the homomorphisms $\phi,\psi:C(U_q(2)) \to C(\bbbt)$ defined by
\begin{align*}
\phi &: a\mapsto 1\quad,\quad b\mapsto 0\quad,\quad D\mapsto\bbt\,;\\
\psi &: a\mapsto\bbt\quad,\quad b\mapsto 0\quad,\quad D\mapsto 1\,.
\end{align*}
It is easy to check that $\phi\mbox{ and } \psi$ satisfy the conditions $\Delta_{\bbbt}\circ \phi=(\phi\ots\phi)\circ \Delta$ and $\Delta_{\bbbt}\circ \psi=(\psi\ots\psi)\circ \Delta$ respectively. Then, the homomorphisms $\Phi: C(U_q(2))\to C(\bbbt)\ots C(U_q(2))$ and $\Psi: C(U_q(2))\to C(\bbbt)\ots C(U_q(2))$ defined by the following formulae
\[
\Phi(x)=(\phi\ots \mbox{id})\Delta\quad,\quad\Psi(x)=(\psi\ots \mbox{id})\Delta
\]
are two different actions of $\bbbt$ on the quantum group $U_q(2)$. The corresponding homogeneous spaces $U_q(2)/_{\phi}\bbbt$ and $U_q(2)/_{\psi}\bbbt$ are defined as follows
\begin{align*}
C(U_q(2)/_{\phi}\bbbt) &:= \{x\in C(U_q(2))\,: (\phi\ots \mbox{id})\Delta(x)=1\ots x\},\\
C(U_q(2)/_{\psi}\bbbt) &:= \{x\in C(U_q(2))\,: (\psi\ots \mbox{id})\Delta(x)=1\ots x\}.
\end{align*}
The maps $E_{\phi}:C(U_q(2)) \xrightarrow{} C(U_q(2)/_{\phi}\bbbt)$ and $E_{\psi}:C(U_q(2)) \xrightarrow{} C(U_q(2)/_{\psi}\mathbb{T})$ defined by the following
\begin{IEEEeqnarray}{lcl}\label{expectation}
E_{\phi}(x)=((h_{\bbbt}\circ\phi)\ots \mbox{id})\circ\Delta\quad,\quad E_{\psi}(x)=((h_{\bbbt}\circ\psi)\ots \mbox{id})\circ \Delta
\end{IEEEeqnarray}
are the conditional expectations for the unital $C^*$-inclusions $\,C(U_q(2)/_{\phi}\bbbt)\subset C(U_q(2))$ and $C(U_q(2)/_{\psi}\mathbb{T})\subset C(U_q(2))$ respectively.
 
\bppsn\label{index1}
The conditional expectation $E_{\phi}:C(U_q(2))\rightarrow C(U_q(2)/_{\phi}\bbbt)$ has infinite Watatani index.
\eppsn
\prf 
Consider the sequence $x_n:=\sum_{j=0}^{n-1}D^j$ in $C(U_q(2))$. If we assume that $E_\phi$ is of index-finite type, then by \Cref{indexcondition} there exists $0<c\le 1$ such that $E(x_n^*x_n)\ge c(x_n^*x_n)$ for all $n$. That is, 
\[
E(x_n^*x_n) \ge c\Big(n+\sum_{j=1}^{n-1}(n-j)(D^j+(D^*)^j)\Big).
\]
Using \Cref{expectation}, we have $E_\phi(D^k)=E_\phi((D^*)^k)=0$ for all $1\leq k\leq n-1$. Therefore, if we choose $\,\xi_n=\sum_{i=0}^{n-1}e_{0,0,i}\in \ell^2(\bbn\times\bbz\times\bbz)$ for $n\in\bbn$, then we get the following
\begin{IEEEeqnarray*}{rcl}
\frac{n}{c}-n\,\,\,&\geq& \,\,\Big|\Big|\sum_{j=1}^{n-1}(n-j)(D^j+(D^*)^j)\,\Big|\Big|\\
&\geq& \,\,\Big\langle \pi\Big(\sum_{j=1}^{n-1}(n-j)(D^j+(D^*)^j)\Big)\xi_n, \xi_n\Big\rangle\\
&=& \,\,\sum_{j=1}^{n-1}(n-j)\Big\langle \pi(D)^{j}\xi_n+\pi(D^*)^j\xi_n, \xi_n\Big\rangle\\
&=& \,\,\sum_{j=1}^{n-1}\frac{(n-j)}{n}\Big\langle \sum_{i=0}^{n-1}e_{0,0,i+j}\,+\sum_{i=0}^{n-1}e_{0,0,i-j}\, , \sum_{i=0}^{n-1}e_{0,0,i}\Big\rangle \\
&=& \,\,\sum_{j=1}^{n-1}\frac{(n-j)}{n}\Big(\Big\langle \sum_{i=j}^{n+j-1}e_{0,0,i}+\sum_{i=-j}^{n-j-1}e_{0,0,i}\, , \sum_{i=0}^{n-1}e_{0,0,i}\Big\rangle\Big) \\
&=& \,\,\sum_{1<2j\le n-1,\,j\in\bbn}\frac{(n-j)}{n}\Big\langle\Big(\sum_{i=0}^{j-1}e_{0,0,i}+2\sum_{i=j}^{n-j-1}e_{0,0,i}+\sum_{i=n-j}^{n-1}e_{0,0,i}\Big)\,,\sum_{i=0}^{n-1}e_{0,0,i}\Big\rangle\\
& & \quad+\sum_{n-1<2j\le2n-2,\,j\in\bbn}\frac{(n-j)}{n}\Big\langle\Big(\sum_{i=0}^{n-j-1}e_{0,0,i}+\sum_{i=j}^{n-1}e_{0,0,i}\Big)\,,\sum_{i=0}^{n-1}e_{0,0,i}\Big\rangle\\
&=& \,\,\frac{2}{n}\sum_{j=1}^{n-1}(n-j)^2\,.
\end{IEEEeqnarray*}
Thus, we get that $\frac{n}{c}\ge n+\frac{2n^2-3n+1}{3}$, and consequently $c\le\frac{3n}{2n^2+1}$ for all $n$. This contradicts our assumption that $c>0$ which says that no non-zero $c(E_\phi)$ in \Cref{indexcondition} can exist, and hence $E_\phi$ is not of index-finite type.\qed

\bppsn
The conditional expectation $E_{\psi}:C(U_q(2))\rightarrow C(U_q(2)/_{\psi}\mathbb{T})$ has infinite Watatani index. 
\eppsn
\prf Consider the sequence $\,y_n:=\sum_{j=0}^{n}(b)^j$ in $C(U_q(2))$. If we assume that $E_\psi$ is of index-finite type, then by \Cref{indexcondition} there is a constant $0<c\leq 1$ such that $E_\psi(y_n^*y_n)\ge c(y_n^*y_n)$ for all $n\in\bbn$. Using \Cref{expectation}, we have $E_\psi(b^l)=E_\psi((b^*)^l)=0$ for all $1\leq l\leq n$. Thus we have
\[
\frac{1}{c}\sum_{j=0}^{n}(b^*b)^j\,\,\ge\,\,\sum_{j=0}^{n}(b^*b)^j+  \sum_{i=1}^{n} \Big(\sum_{j=0}^{n-i}(b^*b)^j\Big)\Big((b^*)^i+b^i\Big)
\]
Therefore, if we choose $\xi_n =\sum_{k=0}^{n}e_{0,k,0} \in \ell^2(\bbn\times\bbz\times\bbz)$ for $n\in \bbn$, then we have the following
\begin{IEEEeqnarray*}{rcl}
& & \frac{1}{c}\Big|\Big|\sum_{j=0}^{n}(b^*b)^j\Big|\Big|\\
&\geq& \,\,\Big\langle\pi\Big( \sum_{j=0}^{n}(b^*b)^j+  \sum_{i=1}^{n} \Big(\sum_{j=0}^{n-i}(b^*b)^j\Big)\Big((b^*)^i+b^i\Big)\xi_n\,,\,\xi_n\Big\rangle\\
&=& \,\,\Big\langle\sum_{j=0}^{n}\pi(b^*b)^j \xi_n\,,\,\xi_n\Big\rangle\,+ \Big\langle \sum_{i=1}^{n} \Big(\sum_{j=0}^{n-i}\pi(b^*b)^j\Big)\Big(\pi(b^*)^i+\pi (b)^i\Big) \xi_n\,,\,\xi_n\Big\rangle\\
&=& \,\,(n+1)+\Big\langle \sum_{i=1}^{n} \frac{n-i+1}{n+1}\Big(\sum_{k=0}^{n}e_{0,k-i,0}+\sum_{k=0}^{n}e_{0,k+i,0} \Big)\,,\,\sum_{k=0}^{n}e_{0,k,0} \Big\rangle\\
&=& \,\,(n+1)+\Big\langle \sum_{i=1}^{n}\frac{n-i+1}{n+1} \Big(\sum_{k=-i}^{n-i}e_{0,k,0}+\sum_{k=i}^{n+i}e_{0,k,0} \Big)\,,\,\sum_{k=0}^{n}e_{0,k,0} \Big\rangle\\
&=& \,\,(n+1)+\Big\langle \sum_{1 \le i \le n-i,i \in \bbn}\frac{n-i+1}{n+1} \Big(\sum_{k=0}^{i-1}e_{0,k,0}+2\sum_{k=i}^{n-i}e_{0,k,0}+\sum_{k=n-i+1}^{n}e_{0,k,0} \Big)\,,\,\sum_{k=0}^{n}e_{0,k,0} \Big\rangle\\
&=& \,\,\Big\langle \sum_{n-i<i\le n,n\in \bbn}\frac{n-i+1}{n+1} \Big(\sum_{k=0}^{n-i}e_{0,k,0}+\sum_{k=i}^{n}e_{0,k,0} \Big)\,,\,\sum_{k=0}^{n}e_{0,k,0} \Big\rangle\\
&=& \,\,(n+1)+\frac{2}{n+1}(\sum_{i=1}^n(n-i+1)^2)\,.
\end{IEEEeqnarray*}
This gives us that $(n+1)+\frac{2}{n+1}(\sum_{i=1}^n(n-i+1)^2)\leq \frac{n+1}{c}$ since $||b^*b||\le 1$, and consequently $c\le \frac{3n+3}{2n^2+4n+3}$ for all $n$. This contradicts our assumption that $c>0$ which says that no non-zero $c(E_\psi)$ in \Cref{indexcondition} can exist, and therefore $E_\psi$ is not of index-finite type.\qed


\newsection{\texorpdfstring{$\mathbb{S}^1\times\mathbb{S}^1\times\mathbb{S}^1\,$}{}-equivariant Dirac operator}\label{Sec5}

In this section, we explore a noncommutative geometric aspect of the compact quantum group $U_q(2)$ along the line of \cite{Con1, Con2}.
\begin{dfn}
Let $\mathscr{A}$ be a unital $C^*$-algebra and $\cla\subseteq\mathscr{A}$ be a unital dense $\star$-subalgebra. An odd spectral triple is a tuple $(\cla,\clh,\mathcal{D})$, where $\clh$ is a separable Hilbert space on which $\cla$ acts as bounded operators, $\mathcal{D}$ is an unbounded self-adjoint operator with compact resolvent such that the commutator $[\mathcal{D},a]$ extends to a bounded operator on $\clh$ for all $a\in\cla$.
\end{dfn}

If there is a $\bbz_2$-grading operator $\gamma$ acting on $\clh$ such that it commutes with $a$ for all $a\in\cla$ and anticommutes with $\mathcal{D}$, then $(\cla,\clh,\mathcal{D},\gamma)$ is called an {\it even} spectral triple. The operator $\mathcal{D}$ is widely referred as the Dirac operator. Since the kernel of $\mathcal{D}$ is finite dimensional, without loss of generality, one can assume that $\mathcal{D}$ has trivial kernel (see Page $316$ in \cite{Con1} or Page $446$ in \cite{GVF}). If for $p>0,\,|\mathcal{D}|^{-p}$ lies in the Dixmier ideal $\mathcal{L}^{(1,\infty)}\subseteq\mathcal{B}(\clh)$, then we say that the spectral triple is $p^+$-summable \cite{Con1}.

An even spectral triple $(\cla,\clh,\mathcal{D},\gamma)$ induces a $\mathcal{K}$-homology class $[(\cla,\clh,F,\gamma)]$, where $F=\mathcal{D}|\mathcal{D}|^{-1}$, in $K^0(\mathscr{A})$ consisting of even Fredholm modules. Given a projection $P\in M_n(\mathscr{A})$, define $\clh_n=\clh\otimes \bbc^n,\,\gamma_n=\gamma\otimes I_n,\,F_n=F\otimes I_n,\,P^+=\frac{1+\gamma_n}{2}P$, and $P^-=\frac{1-\gamma_n}{2}P$. Then, $\clh_n$ decomposes as $\clh_n^+\oplus\clh_n^-$ under the grading operator $\gamma_n$ and the operator $P^-F_nP^+:P^+\clh_n^+\longrightarrow P^-\clh_n^-$ is a Fredholm operator. Index of this Fredholm operator is the value of the $K_0-K^0$ pairing $\langle[P],[(\cla,\clh,F,\gamma)]\rangle$. An even spectral triple $(\cla,\clh,\mathcal{D},\gamma)$ is called nontrivial if $\langle[P],[(\cla,\clh,F,\gamma)]\rangle$ is nonzero for some $[P]\in K_0(\mathscr{A})$. The requirement of nontrivial pairing is very crucial and see \cite{CP}, for instance, for the reason.

\begin{dfn}
Let $\mathbb{G}$ be a compact group acting on a unital $C^*$-algebra $\mathscr{A}$ via the strongly continuous action $\alpha$ so that we have a $C^*$-dynamical system $(\mathscr{A},\mathbb{G},\alpha)$. A covariant representation $(\pi,\mathbb{U})$ of $(\mathscr{A},\mathbb{G},\alpha)$ consists of a unital $\star$-representation $\pi:\mathscr{A}\longrightarrow\mathbb{B}(\mathcal{H})$ and a unitary representation $\mathbb{U}:\mathbb{G}\longrightarrow\mathcal{U}(\mathcal{H})$ such that $\pi(\alpha_g(x))=\mathbb{U}_g\pi(x)\mathbb{U}_g^*$ for all $x\in\mathscr{A}$ and $g\in\mathbb{G}$.
\end{dfn}

\bdfn\label{definition}
Let $(\mathscr{A},\mathbb{G},\alpha)$ be a $C^*$-dynamical system. An even $\mathbb{G}$-equivariant spectral triple for $(\mathscr{A},\mathbb{G},\alpha)$ is a tuple $(\pi,\mathbb{U},\mathscr{H},\mathscr{D},\gamma)$, where
\begin{enumerate}[$(i)$]
\item $(\pi,\mathbb{U})$ is a covariant representation of $(\mathscr{A},\mathbb{G},\alpha)$ on the separable Hilbert space $\mathscr{H}$,
\item $\pi$ is faithful,
\item $(\pi,\mathscr{H},\mathscr{D},\gamma)$ is an even spectral triple for $\mathscr{A}$,
\item $\mathbb{U}_g\mathscr{D}\mathbb{U}_g^*=\mathscr{D}$ and $[\mathbb{U}_g,\gamma]=0$ for all $g\in\mathbb{G}$.
\end{enumerate}
\edfn

Detail discussion on equivariant spectral triple can be found in \cite{Sitarz}. If all the notations are clear from the context, then we simply say that the Dirac operator $\mathscr{D}$ is $\mathbb{G}$-equivariant. If a compact group acts on a unital $C^*$-algebra, one demands the condition of ``equivariance'' in order to capture the geometry of the underlying dynamical system. In \cite{GS-20212aa}, a $4^+$-summable even spectral triple has been constructed on $U_q(2)$ that is equivariant under its own comultiplication action. Our purpose here is to construct a torus-equivariant Dirac operator on $U_q(2)$.
\smallskip

Let $\clh:=\ell^2(\bbn)\otimes\ell^2(\bbz)\otimes\ell^2(\bbz)$ and recall the faithful representation $\pi:C(U_q(2))\rightarrow\clb(\clh)$ described in \Cref{representation}. The compact Lie group $\mathbb{G}:=\mathbb{T}^3$ acts on $C(U_q(2))$ by the strongly continuous homomorphism $\alpha:\mathbb{T}^3\longrightarrow\mbox{Aut}(C(U_q(2)))$ defined by
\begin{IEEEeqnarray}{rCl}\label{the action}
a\longmapsto z_1a\quad,\quad b\longmapsto z_2b\quad,\quad D\longmapsto z_3D
\end{IEEEeqnarray}
for $(z_1,z_2,z_3)\in\mathbb{T}^3$. We have a unitary representation $\mathbb{U}$ of $\mathbb{T}^3$ on $\clh$ given by $\mathbb{U}_{(z_1,z_2,z_3)}=z_1^N\otimes z_2^N\otimes z_3^N$, where $N$ denotes the number operator $e_n\longmapsto ne_n$ on $\ell^2(\bbn)$ or $\ell^2(\bbz)$. It is easy to check that $\pi(\alpha_{(z_1,z_2,z_3)}(x))=\mathbb{U}_{(z_1,z_2,z_3)}\pi(x)\mathbb{U}_{(z_1,z_2,z_3)}^*$ for all $x\in C(U_q(2))$ and $(z_1,z_2,z_3)\in\mathbb{T}^3$. That is, $(\mathbb{U},\pi,\clh)$ is a covariant representation of the $C^*$-dynamical system $(C(U_q(2)),\mathbb{T}^3,\alpha)$.

Let $\cla_f$ be the $\star$-subalgebra of $C(U_q(2))$ generated by $a,b,D$. In order to avoid any possible confusion with inner-product, in this particular section we denote the basis element $\langle n,m,r,s\rangle$ (see \Cref{basis}) by $a_nb^m(b^*)^rD^s$, that is,
\begin{IEEEeqnarray}{rCl}\label{the basis}
a_nb^m(b^*)^rD^s:=\begin{cases}
                            a^n b^m(b^*)^rD^s & \mbox{ if } n\geq 0\,,\cr
                            (a^*)^{-n} b^m(b^*)^rD^s & \mbox{ if } n< 0\,.\cr
                           \end{cases}
\end{IEEEeqnarray} 
We are going to investigate a $\mathbb{T}^3$-equivariant, $K$-homologically non-trivial, finitely summable, even spectral triple on $\cla_f$.

Consider the unbounded operator $R:e_{i,j}\longmapsto d(i,j)e_{i,j}$ acting on $\ell^2(\bbn)\otimes\ell^2(\bbz)$ and defined by the following\,:
\begin{IEEEeqnarray}{lCl}\label{the first sequence}
d{(i,j)}=\begin{cases}
i+j& \mbox{ if }\,\,j\geq 0,\cr
-i+j& \mbox{ if }\,\,j<0,\cr
\end{cases}
\end{IEEEeqnarray}
Now, consider the unbounded operator $T:e_{i,j,k}\longmapsto d(i,j,k)e_{i,j,k}$ on $\mathcal{H}$ given by $d(i,j,k)=d(i,j)+\sqrt{-1}k$. Hence, $T=R\otimes 1+\sqrt{-1}(1\otimes 1\otimes N)$. Thus, from \Cref{the first sequence} we get that
\begin{IEEEeqnarray}{lCl}\label{the second sequence}
d(i,j,k)=\begin{cases}
i+j+\sqrt{-1}k& \mbox{ if }\,\,j\geq 0,\cr
-i+j+\sqrt{-1}k& \mbox{ if }\,\,j<0,\cr
\end{cases}
\end{IEEEeqnarray}

We first deal with the case of $\theta\neq 0,1$, where $\theta=\frac{1}{\pi}\arg{(q)}$.

\blmma\label{bounded commutator}
For each $x\in\mathcal{A}_f$, the operators $[T,x]$ and $[T^*,x]$ extend to bounded operators on $\clh$.
\elmma
\prf We only show that $[T,x]$, where $x\in\{a,a^*,b,b^*,D,D^*\}\subseteq\mathcal{A}_f$, extends to a bounded operator on $\clh$. From \Cref{the second sequence}, we get the following set of six equality\,:
\begin{align*}
[T,D](e_{i,j,k}) &= \big(d(i,j,k+1)-d(i,j,k)\big)e^{-2\pi\theta j\sqrt{-1}}e_{i,j,k+1}\,\,,\\
[T,D^*](e_{i,j,k}) &= \big(d(i,j,k-1)-d(i,j,k)\big)e^{2\pi\theta j\sqrt{-1}}e_{i,j,k-1}\,\,,\\
[T,a](e_{i,j,k}) &= \big(d(i+1,j,k)-d(i,j,k)\big)\sqrt{1-|q|^{2i+2}}e_{i+1,j,k}\,\,,\\
[T,a^*](e_{i,j,k}) &= \big(d(i-1,j,k)-d(i,j,k)\big)\sqrt{1-|q|^{2i}}e_{i-1,j,k}\,\,,\\
[T,b](e_{i,j,k}) &= \big(d(i,j+1,k)-d(i,j,k)\big)q^ie_{i,j+1,k}\,\,,\\
[T,b^*](e_{i,j,k}) &= \big(d(i,j-1,k)-d(i,j,k)\big)\overline{q}^ie_{i,j-1,k}\,\,.
\end{align*}
Thus, all the commutators are weighted shift either on $\ell^2(\bbn)$ or $\ell^2(\bbz)$. Observe that
\[
|q|^i \,|d(i,j+1,k)-d(i,j,k)|=\begin{cases}
|q|^i& \mbox{ for }\,\,j\neq -1\,;\\
|2i+1|\,|q|^i& \mbox{ for }\,\,j=-1\,.
\end{cases}
\]
and
\[
|\overline{q}|^i \,|d(i,j-1,k)-d(i,j,k)|=\begin{cases}
|q|^i& \mbox{ for }\,\,j\neq 0\,;\\
|2i+1|\,|q|^i& \mbox{ for }\,\,j=0\,.
\end{cases}
\]
This gives the boundedness of the commutators $[T,b]$ and $[T,b^*]$. In all the other cases, the respective weight functions are bounded by $1$, and this completes the proof.\qed

Define the faithful representation $\pi_{eq}$ of $\,C(U_q(2))$ on $\clh\otimes\mathbb{C}^2$ by
\[
\pi_{eq}(x)=\left[ {\begin{matrix}
   \pi(x)  & 0\\
  0 & \pi(x) \\
  \end{matrix} } \right]\,\,.
\]
Let 
\[
\mathscr{D}=\left[ {\begin{matrix}
   0  & T^*\\
  T & 0 \\ 
  \end{matrix} } \right]\quad \mbox{ and }\quad \gamma= \left[ {\begin{matrix}
   1 & 0\\
  0 & -1 \\
  \end{matrix} } \right].
\]
Immediately from \Cref{bounded commutator}, it follows that $[\mathscr{D},\pi_{eq}(x)]$ extends to a bounded operator on $\clh\otimes\bbc^2$ for each $x\in\mathcal{A}_f$. The operator $\mathscr{D}$ becomes essentially self-adjoint. The unique self-adjoint extension of $\mathscr{D}$ is again denoted by the same symbol, as is customary in noncommutative geometry.

\blmma\label{compact resolvent}
The unbounded operator $\mathscr{D}$ has compact resolvent.
\elmma
\prf Let $\mathcal{W}=\mathrm{span}\{(e_{0,0,0},0),(0,e_{0,0,0})\}$ be the two dimensional subspace of $\clh\oplus\clh$ which is the kernel of $\mathscr{D}$. On $\mathcal{W}^\perp$ in $\clh\otimes\mathbb{C}^2$, observe the following,
\[|T^*|^{-1}(e_{i,j,k})=|T|^{-1}(e_{i,j,k})=\frac{1}{\sqrt{(i+|j|)^2+k^2}}\,e_{i,j,k}\,.\] Thus, we have $|\mathscr{D}|^{-1}=|T|^{-1}\otimes I_2$ on $\mathcal{W}^\perp$, and hence is a compact operator on $\clh\otimes\bbc^2$.\qed

Therefore, $\mathscr{D}$ is a Dirac operator on $\mathcal{A}_f$. We now concentrate on its summability.

\bppsn\label{summability}
One has $|\mathscr{D}|^{-3}\in\mathcal{L}^{(1,\infty)}$, the ideal of the Dixmier traceable operators.
\eppsn
\prf Observe that $\mathscr{D}=R\otimes 1\otimes\sigma_1+1\otimes N\otimes\sigma_2$, where the matrices $\sigma_1=\begin{pmatrix}
0 & 1\\
1 & 0
\end{pmatrix}$ and $\sigma_2=\begin{pmatrix}
0 & -i\\
i & 0
\end{pmatrix}$ are the Pauli spin matrices. The grading operator $\gamma$ is the third Pauli spin matrix $\sigma_3$. Therefore, $\mathscr{D}$ acting on $\clh\otimes\bbc^2$ is of the form $\mathscr{D}^\prime\otimes 1\otimes\sigma_1+1\otimes N\otimes\sigma_2$, where $\mathscr{D}^\prime$ is exactly the $\mathbb{T}^2$-equivariant Dirac operator for $SU_q(2)$ acting on $\ell^2(\bbn)\otimes\ell^2(\bbz)$ obtained in \cite{CP1}, and $N$ is the standard Dirac operator for $\mathbb{T}$ acting on $\ell^2(\bbz)$. It is known that $\mathscr{D}^\prime$ is $2^+$-summable by \cite{CP1}, i,e. $Tr_\omega(|\mathscr{D}^\prime|^{-2})<\infty$ ($Tr_\omega$ denotes the Dixmier trace), and $N$ is $1^+$-summable with $Tr_\omega(|N|^{-1})=1$. Therefore, by the result of Connes (Page $576$ in \cite{Con1}), we get that
\begin{IEEEeqnarray*}{lCl}
Tr_\omega(|\mathscr{D}|^{-3}) &=& \frac{\Gamma(2)\Gamma\left(\frac{3}{2}\right)}{\Gamma\left(\frac{5}{2}\right)}Tr_\omega(|\mathscr{D}^\prime|^{-2})Tr_\omega(|N|^{-1})\\
&=& \frac{2}{3}\,Tr_\omega(|\mathscr{D}^\prime|^{-2})\,,
\end{IEEEeqnarray*}
which completes the proof.\qed

\begin{rmrk}\rm
When $\theta=\frac{1}{\pi}\arg{(q)}$ and $\theta\notin\{0,1\}$, it is not true that as a $C^*$-algebra one has $C(U_q(2))=C(SU_q(2))\otimes C(\mathbb{T})$. The tensor product `$\otimes$' must be replaced by the braided tensor product `$\boxtimes$'. Therefore, the realization of the Dirac operator $\mathscr{D}$ in terms of $\mathscr{D}^\prime$ is not automatic {\em a priori}.
\end{rmrk}

It follows that $\,\mathbb{U}_{(z_1,z_2,z_3)}T\mathbb{U}_{(z_1,z_2,z_3)}^*=T$ and $\,\mathbb{U}_{(z_1,z_2,z_3)}T^*\mathbb{U}_{(z_1,z_2,z_3)}^*=T^*$. Hence, $\mathscr{D}$ is equivariant under the $\mathbb{T}^3$-action defined in \Cref{the action}. Combining this with Lemma (\ref{bounded commutator}, \ref{compact resolvent}) and \Cref{summability}, we see that all the requirements in \Cref{definition} are now satisfied to conclude the following theorem.

\bthm\label{spectral triple}
The tuple $(\mathcal{A}_f\,,\clh\otimes\mathbb{C}^2\,,\pi_{eq}\,,\mathscr{D}\,,\gamma)$ is a $\mathbb{T}^3$-equivariant, $3^+$-summable, even spectral triple for the compact quantum group $U_q(2)$ for $|q|\neq 1$. 
\ethm

It is the next theorem, the irrationality of $\theta$ again plays a crucial role. This is because the $K$-groups in \cite{GS-20212aa} have been computed under this additional hypothesis, and the reason is its connection with the noncommutative torus $\mathbb{A}_\theta$.

\bthm\label{Chern}
The Chern character of the spectral triple constructed in \Cref{spectral triple} is non-trivial when $\theta$ is irrational.
\ethm
\prf We calculate the index of the Fredholm operator
\[
P_\theta T|T|^{-1}P_\theta:P_\theta\clh\longrightarrow P_\theta\clh
\]
where $P_\theta=[p\otimes p_\theta]$, a generator of $K^0(C(U_q(2)))$. Here, $p:=|e_0\rangle\langle e_0|$ is the rank one projection acting on $\ell^2(\bbn)$ and $p_\theta$ is the Powers-Rieffel projection in the noncommutative torus $\mathbb{A}_\theta$ with trace $\theta$ (see Theorem $3.5$ in \cite{GS-20212aa} for detail). Let us first investigate how the following operator
\[
(p\otimes 1\otimes 1)T|T|^{-1}(p\otimes 1\otimes 1):\clh\longrightarrow\clh
\]
acts on $\clh$. For the standard orthonormal basis $\{e_{i,j,k}\}$ of $\clh$ we get that
\begin{IEEEeqnarray*}{lCl}
(p\otimes 1\otimes 1)T|T|^{-1}(p\otimes 1\otimes 1)(e_{i,j,k}) &=& (p\otimes 1\otimes 1)T|T|^{-1}(e_{0,j,k})\\
&=& (p\otimes 1\otimes 1)\Big(\frac{d(0,j,k)}{\sqrt{j^2+k^2}}\,e_{0,j,k}\Big)\\
&=& \frac{j+\sqrt{-1}k}{\sqrt{j^2+k^2}}\,e_{0,j,k}\,.
\end{IEEEeqnarray*}
Therefore, $(p\otimes 1\otimes 1)T|T|^{-1}(p\otimes 1\otimes 1)$ is the operator $p\otimes\widetilde{T}$ acting on $\clh$, where $\widetilde{T}:\ell^2(\bbz^2)\longrightarrow\ell^2(\bbz^2)$ is given by $e_{j,k}\longmapsto \frac{j+\sqrt{-1}k}{\sqrt{j^2+k^2}}e_{j,k}$. Thus, $\widetilde{T}=T_{m+\sqrt{-1}n}|T_{m+\sqrt{-1}n}|^{-1}$, where $T_{m+\sqrt{-1}n}(e_{m,n})=(m+\sqrt{-1}n)e_{m,n}$ acting on $\ell^2(\bbz^2)$. It is known that the Fredholm operator
\[
p_\theta\widetilde{T}p_\theta:p_\theta\ell^2(\bbz^2)\longrightarrow p_\theta\ell^2(\bbz^2)
\]
associated with the Dirac operator $\begin{pmatrix}
0 & T_{m-\sqrt{-1}n}\\
T_{m+\sqrt{-1}n} & 0
\end{pmatrix}$ for the noncommutative torus has nontrivial index \cite{Had}. Observe that
\begin{IEEEeqnarray*}{lCl}
P_\theta T|T|^{-1}P_\theta &=& (1\otimes p_\theta)(p\otimes\widetilde{T})(1\otimes p_\theta)\\
&=& p\otimes\big(p_\theta\widetilde{T}p_\theta\big)
\end{IEEEeqnarray*}
as an operator acting on $\,p\ell^2(\bbn)\otimes p_\theta\ell^2(\bbz^2)$. Since $p$ is a rank one projection acting on $\ell^2(\bbn)$, we conclude that the index of $P_\theta T|T|^{-1}P_\theta$ is nonzero.\qed

We finally comment for the cases $\theta=0,1$. In these cases, $q$ becomes a real deformation parameter, and by Theorem $2.1$ in \cite{Z} we know that $C(U_q(2))\cong C(SU_q(2))\otimes C(\mathbb{T})$ as $C^*$-algebras. Therefore, Theorem \ref{spectral triple} obviously follows. Moreover, in this case $[p\otimes Bott]$ is a generator of $K^0(C(U_q(2)))$, and Theorem \ref{Chern} also holds if one replaces $p_\theta$ by ``$Bott$'' (see section $3$ in \cite{GS-20212aa} for detail).
\medskip

Recall that a spectral triple such that $[\mathscr{D},\pi(x)]=0$ implies $\,x\in\bbc$ is called non-degenerate. The above spectral triple is not non-degenerate as one can quickly observe that $[\mathscr{D},bb^*]$ is zero by \Cref{representation,the second sequence}. In fact, no $\mathbb{T}^3$-equivariant Dirac operator can be non-degenerate for the representation described in \Cref{representation}. Same is the situation for the $\mathbb{T}^2$-equivariant Dirac operator on $SU_q(2)$ \cite{CP1}. We conclude this section by identifying the kernel of the unbounded derivation $d:=[\mathscr{D},\,\,]$ in $\mathcal{A}_f$.

\bppsn
The kernel of the unbounded derivation $\,d:=[\mathscr{D},\,\,]$ in $\mathcal{A}_f$ is the subspace $\,\{x=\sum_{j=0}^n\gamma_jb^j(b^*)^j:\gamma_j\in\bbc,\,n\in\bbn\}$ of $\mathcal{A}_f$.
\eppsn
\prf In view of \Cref{the basis}, consider any finite subset $F$ of $\bbz\times\bbn\times\bbn\times\bbz$ and scalars $\,C_{n,m,r,s}\neq 0$. Take $\,x=\sum_{(n,m,r,s)\in F}C_{n,m,r,s}\,a_nb^m(b^*)^rD^s\in\mathcal{A}_f$. Consider the commutator $[T,x]$. By \Cref{representation,the second sequence}, we get the following
\begin{IEEEeqnarray*}{lCl}
&  & [T,x](e_{i,j,k})\\
&=& \sum_{(n,m,r,s)\in F}C_{n,m,r,s}\,e^{-2js\pi\sqrt{-1}\theta}(\overline{q})^{ir}q^{im}\alpha(n)\,\big(d(i+n,j-r+m,k+s)-d(i,j,k)\big)e_{i+n,j-r+m,k+s}
\end{IEEEeqnarray*}
where we have
\begin{IEEEeqnarray}{lCl}\label{alpha}
\alpha(n)=\begin{cases}
\prod_{\xi=1}^n\sqrt{1-|q|^{2(i+\xi)}} & \mbox{ if }\,\,n\geq 1\,,\cr
1 & \mbox{ if }\,\,n=0\,,\cr
\prod_{\xi=0}^{-n-1}\sqrt{1-|q|^{2(i-\xi)}} & \mbox{ if }\,\,n<0\,.\cr
\end{cases}
\end{IEEEeqnarray}
Let $\eta=\max\{|n|+|s|+m+r\}\in\bbn$ and let $(n^\prime,m^\prime,r^\prime,s^\prime)\in F$ be such that $|n^\prime|+|s^\prime|+m^\prime+r^\prime=\eta$. Fix this tuple in $F$. Let $F^\prime=\{(n,m,r,s)\in F:|n|+|s|+m+r=\eta\}$. Then, $F^\prime\neq\emptyset$ as $(n^\prime,m^\prime,r^\prime,s^\prime)\in F^\prime$. Let $F^{\prime\prime}=\{(m,r):(n^\prime,m,r,s^\prime)\in F^\prime\}$ and take $\xi=\max\{m-r:(m,r)\in F^{\prime\prime}\}$. Now take $F^{\prime\prime\prime}=\{(m,r)\in F^{\prime\prime}:\xi=m-r\}$. We claim that $F^{\prime\prime\prime}$ is singleton.
Suppose there exist $(m_1,r_1)\neq(m_2,r_2)$ in $F^{\prime\prime\prime}$. Then,
\begin{IEEEeqnarray*}{lCl}
\xi=m_1-r_1=m_2-r_2\,,\\
\eta-|n^\prime|=m_1+r_1=m_2+r_2\,,
\end{IEEEeqnarray*}
with $\xi,\eta-|n^\prime|$ fixed integer. This proves that $\,m_1=m_2$ and $r_1=r_2$. Thus, the set $F^{\prime\prime\prime}=\{(m_0,r_0)\}$ is singleton. We have the following cases.
\smallskip

\noindent\textbf{Case $1$~:} Suppose that $m_0>r_0$. Observe that for any large enough $i\in\bbn$ with $i+n^\prime\gg 0$,
\begin{IEEEeqnarray*}{lCl}
&  & \langle e_{i+n^\prime,0,s^\prime},[T,x](e_{i,-|m_0-r_0|,0})\rangle\\
&=& \langle e_{i+n^\prime,0,s^\prime},\sum_{(n,m,r,s)\in F}C_{n,m,r,s}e^{2\pi\sqrt{-1}\theta|m_0-r_0|s}(\overline{q})^{ir}q^{im}\alpha(i)\,\big(d(i+n,-|m_0-r_0|+m-r,s)\\
&  & \quad\quad\quad-d(i,-|m_0-r_0|,0)\big)e_{i+n,-|m_0-r_0|+m-r,s}\rangle\\
&=& C_{n^\prime,m_0,r_0,s^\prime}e^{2\pi\sqrt{-1}\theta|m_0-r_0|s}(\overline{q})^{ir_0}q^{im_0}\alpha(i)\big(d(i+n^\prime,0,s^\prime)-d(i,-|m_0-r_0|,0)\big)\,.
\end{IEEEeqnarray*}

\noindent\textbf{Case $2$~:} Suppose that $m_0<r_0$. Observe that for any large enough $i\in\bbn$ with $i+n^\prime\gg 0$,
\begin{IEEEeqnarray*}{lCl}
&  & \langle e_{i+n^\prime,-1,s^\prime},[T,x](e_{i,|m_0-r_0|-1,0})\rangle\\
&=& C_{n^\prime,m_0,r_0,s^\prime}e^{-2\pi\sqrt{-1}\theta(|m_0-r_0|-1)s}(\overline{q})^{ir_0}q^{im_0}\alpha(i)\big(d(i+n^\prime,-1,s^\prime)-d(i,|m_0-r_0|-1,0)\big)\,.
\end{IEEEeqnarray*}
If $s^\prime\neq 0$, we immediately see that the above scalar is non-zero by \Cref{the second sequence,alpha}. For $s^\prime=0$, we observe that
\begin{align*}
d(i+n^\prime,0,s^\prime)-d(i,-|m_0-r_0|,0) &= 2i+n^\prime+m_0-r_0\,,\\
d(i+n^\prime,-1,s^\prime)-d(i,|m_0-r_0|-1,0) &= -2i-n^\prime+m_0-r_0\,.
\end{align*}
Since $n^\prime,\,m_0,\,r_0$ are fixed integers, we can always choose a large enough $i\in\bbn$ to make both the integers $2i+n^\prime+m_0-r_0$ and $2i+n^\prime+m_0-r_0$ nonzero simultaneously. Therefore, the scalars in both the cases are nonzero. Hence, the commutator $[T,x]$ is nonzero in both the above cases and consequently, the derivation $[\mathscr{D},\,\,]$ does not have any kernel. Therefore, the kernel of $[\mathscr{D},\,.\,]$ can be contributed by the only remaining case $m_0=r_0$ and $s^\prime=0$. Now, observe that
\begin{IEEEeqnarray*}{lCl}
\langle e_{\eta+n^\prime,0,0},[T,x](e_{\eta,0,0})\rangle &=& C_{n^\prime,m_0,r_0,0}(\overline{q})^{\eta r_0}q^{\eta m_0}\alpha(\eta)\big(d(\eta+n^\prime,0,0)-d(\eta,0,0)\big)\\
&=& C_{n^\prime,m_0,r_0,0}(\overline{q})^{\eta r_0}q^{\eta m_0}\alpha(\eta)n^\prime\,.
\end{IEEEeqnarray*}
Therefore, if $n^\prime\neq 0$ then $[T,x]$ is again nonzero. Hence, in order to have $[T,x]=0$ we see that $x$ has to be of the form $\,\sum_{j=0}^n\gamma_jb^j(b^*)^j$ with $\gamma_j\in\bbc$ and $n\in\bbn$. In that case, $[T^*,x]$ is also zero, and therefore, $x\in\mbox{Ker}([\mathscr{D},\,.\,])$. The converse is obvious by \Cref{representation,the second sequence}.\qed

\bigskip


\bigskip

\bigskip

\noindent{\sc Debabrata Jana} (\texttt{debabrata.jana05@gmail.com})\\
{\footnotesize Department of Mathematics and Statistics,\\
Indian Institute of Technology, Kanpur,\\
Uttar Pradesh 208016, India}
\end{document}